\documentclass{article}
\usepackage[utf8]{inputenc} 
\usepackage{graphicx} 
\usepackage{booktabs} 
\usepackage{array} 
\usepackage{xcolor}
\usepackage{mparhack}
\usepackage{enumerate}
\usepackage{amsmath}
\usepackage{mathtools}
\usepackage{amssymb}
\usepackage{graphicx}
\usepackage{float}
\usepackage{subfig}
\usepackage{fancyhdr}
\usepackage{siunitx}
\usepackage{dcolumn}
\usepackage{commath}
\usepackage[hidelinks]{hyperref}
\usepackage{tikz}
\usetikzlibrary{arrows}
\usetikzlibrary{matrix}
\usetikzlibrary{decorations.pathreplacing}
\usetikzlibrary{calc}
\newcommand{\KW}{\noindent{\bf Keywords: }}
\newcommand{\AMS}{\medskip\noindent{\bf Mathematics Subject Classification: }}

\makeatletter
\newcommand*\dashline{\rotatebox[origin=c]{90}{$\dabar@\dabar@\dabar@$}}
\makeatother

\begin{document}
\title{Analysis of time series and signals using the Square Wave Method} \author{
Osvaldo Skliar\thanks{Escuela de Inform\'atica, Universidad Nacional, Costa Rica. E-mail: oskliar@una.cr; oskliar@costarricense.cr}
\and Ricardo E. Monge\thanks{Escuela de Ciencias de la Computaci\'on e Inform\'atica, Universidad de Costa Rica, Costa Rica. E-mail: ricardo.mongegapper@ucr.ac.cr}
\and Sherry Gapper\thanks{Universidad Nacional, Costa Rica. E-mail: sherry.gapper.morrow@una.cr}}


\maketitle
\begin{abstract}
The Square Wave Method (SWM), previously introduced for the analysis of signals and images, is presented here as a mathematical tool suitable for the analysis of time series and signals. To show the potential that the SWM has to analyze many different types of time series, the results of the analysis of a time series composed of a sequence of 10,000 numerical values are presented here. These values were generated by using the Mathematical Random Number Generator (MRNG). 
\end{abstract}

\KW  time series analysis; signal analysis; signals with abrupt changes; square wave method; square wave transform

\AMS 62M10, 94A12, 65F99

\newpage
\section{Introduction}
The objective of this article is to address the analysis of time series using a method known as the Square Wave Method (SWM), which was introduced previously for the analysis of signals and images \cite{b0}, \cite{b1}, \cite{b2}.

By using the SWM, consideration is given to the analysis of time series, for which the time interval between any two consecutive values in the series analyzed remains constant. 

If a given time series can be considered to be an adequate representation of a signal during a certain time lapse between an initial instant $t_0$ and a final instant $t_f$, then the analysis of that series, using the SWM, is also an analysis of that signal, within that period of time.

The results obtained after analyzing time series with the SWM are expressed in the frequency domain in a clear, precise and concise way by using a previously introduced mathematical tool, known as the Square Wave Transform (SWT) [1].

Each time series is composed of a sequence of values. One significant characteristic of the SWM is that it takes into account the order in which those values appear. The analysis process of a time series, conducted with the SWM, leads to certain basic components, specific trains of square waves, from which  the time series analyzed can be entirely reconstructed. That is precisely why if a time series may be considered to be a adequate representation of a signal, the analysis of that series using the SWM can also be a valid analysis of that signal.

That significant characteristic of the SWM makes it quite different from other valuable methods, statistical methods in particular, used for the analysis of time series \cite{b3}.

The SWM is reviewed in section 2 to show the relevance of this method as a mathematical tool for the analysis of time series. This information has been included here to provide a solid basis for anyone who has not read previous articles on the SWM.

\section{An example of the analysis of a time series using the SWM}

The procedure to be used to analyze a time series using the SWM is presented in this section.

Suppose that every quarter of a second (\SI{0.25}{\second}) a measurement is made of the difference in the electric potential between two points (point 1 and point 2). To simplify the explanation of the procedure used, suppose that the values measured, expressed in millivolts (\si{\mV}), have been rounded off (for didactic purposes) to whole values. (Later on, the values will no longer be rounded off.) 

Let us admit that the sequence of values measured in millivolts, during \SI{2}{\second}, is the following: 
\begin{equation*}
84, -152, 63, 98, -35, 0, 145, -14
\end{equation*}

Of course, the positive values in this sequence imply that in each case the electric potential of point 1 is greater than that of point 2. Likewise, the negative values in the above sequence imply that in each case the electric potential of point 1 is less than that of point 2. The value of 0 in that sequence of values implies that, under the given precision, the values of the electric potentials of points 1 and 2 are equal.

\begin{figure}[H]
\centering
\includegraphics[width=4.5in]{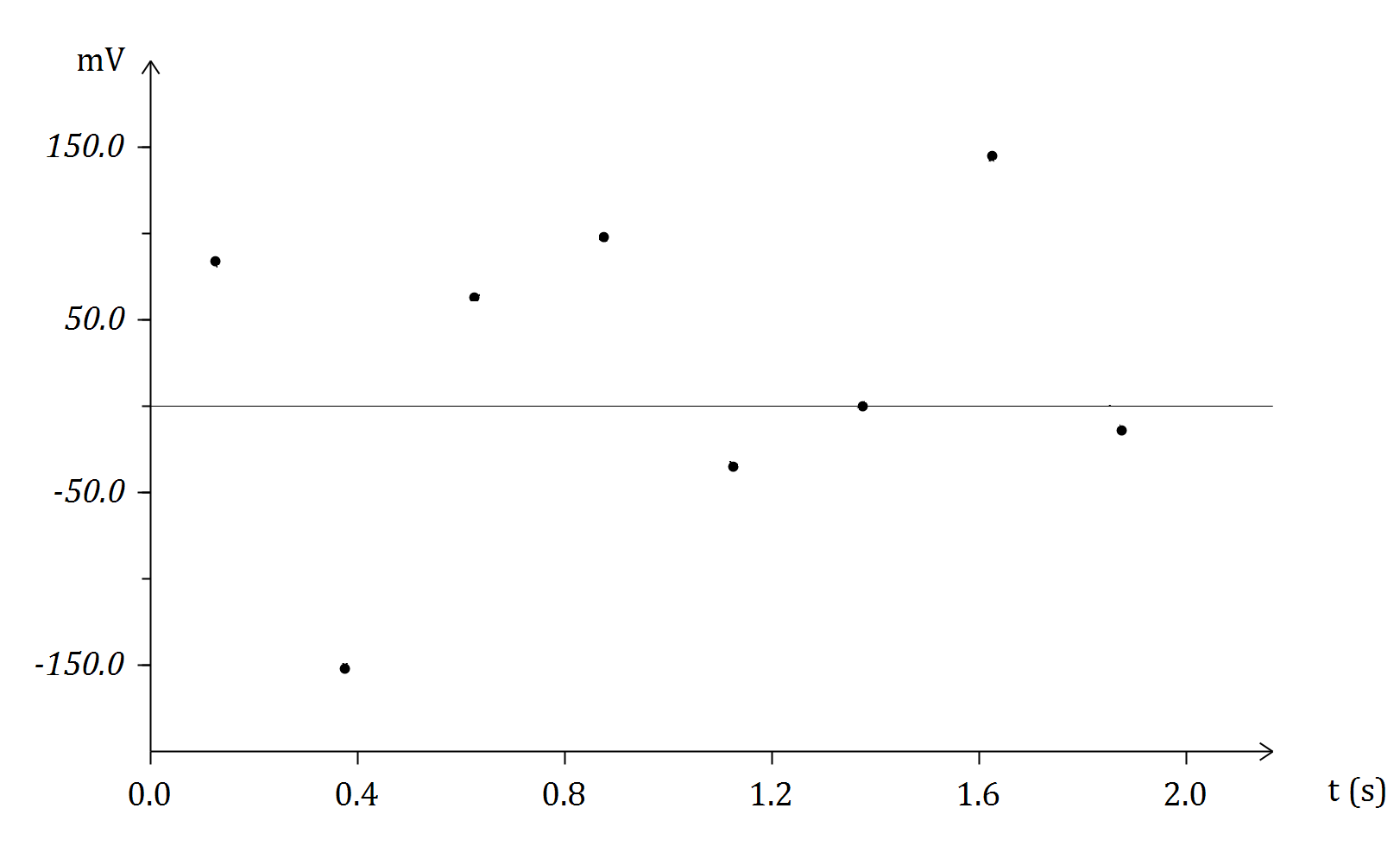}
\caption{Graphic representation of the time series, that is, of the sequence of measured values of differences in electric potential, mentioned above. Note that the time interval between the two consecutive measurements is equal to \SI{0.25}{\second}.}

\label{f1}
\end{figure}

The values of this time series are shown in figure~\ref{f1}.

According to the SWM, each of the measured values in the time series considered can be approximated very precisely by adding up 8 values. Each of these 8 values is located in a corresponding train of square waves.

Reference is made here to the 8 added values and to the 8 corresponding trains of square waves because in this case the time series considered is made up of 8 measured values. If the sequence were composed of 1,000 measured values, it would be necessary to add up the 1,000 values situated in the 1,000 corresponding trains of square waves, to approximate each of those measured values. In general, if that time series is made up of $n$ values, such that $n=1,2,3,\ldots$, then to approximate each of those values, the $n$ values situated in $n$ corresponding trains of square waves must be added up.

An explanation is provided below of how to discover 1) what values must be added up to approximate each measured value; 2) which trains of square waves are being referred to; and 3) where the values that must be summed up are located in those trains of square waves.

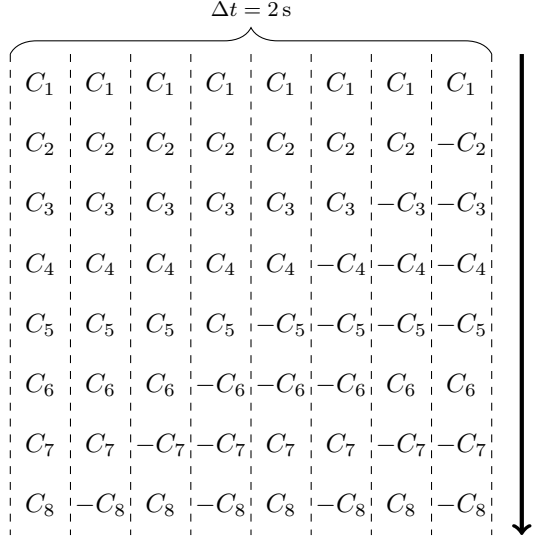
\begin{figure}[H]
\centering
\begin{tikzpicture}
\matrix[inner sep=0pt] (S) 
       [matrix of math nodes,
       nodes={outer sep=0pt,minimum width=8mm,minimum height=8mm}]
      { C_1 & C_1 & C_1 & C_1 & C_1 & C_1 & C_1 & C_1 \\
        C_2 & C_2 & C_2 & C_2 & C_2 & C_2 & C_2 & -C_2 \\
        C_3 & C_3 & C_3 & C_3 & C_3 & C_3 & -C_3 & -C_3 \\
        C_4 & C_4 & C_4 & C_4 & C_4 & -C_4 & -C_4 & -C_4 \\
        C_5 & C_5 & C_5 & C_5 & -C_5 & -C_5 & -C_5 & -C_5 \\
        C_6 & C_6 & C_6 & -C_6 & -C_6 & -C_6 & C_6 & C_6 \\
        C_7 & C_7 & -C_7 & -C_7 & C_7 & C_7 & -C_7 & -C_7 \\
       C_8 & -C_8 & C_8 & -C_8 & C_8 &- C_8 & C_8 & -C_8 \\
      };  
      \draw [decorate,decoration={brace,amplitude=10pt},xshift=0pt,yshift=0pt] 
       (S-1-1.north west) -- (S-1-8.north east)  node [black,midway,yshift=0.6cm]  {\footnotesize $\Delta t=\SI{2}{\second}$};
      \draw[dashed] (S-1-1.north west) --(S-8-1.south west);
      \draw[dashed] (S-1-2.north west) --(S-8-2.south west);
      \draw[dashed] (S-1-3.north west) --(S-8-3.south west);
      \draw[dashed] (S-1-4.north west) --(S-8-4.south west);
      \draw[dashed] (S-1-5.north west) --(S-8-5.south west);
      \draw[dashed] (S-1-6.north west) --(S-8-6.south west);
      \draw[dashed] (S-1-7.north west) --(S-8-7.south west);
      \draw[dashed] (S-1-8.north west) --(S-8-8.south west);
      \draw[dashed] (S-1-8.north east) --(S-8-8.south east);
       \begin{scope}[transform canvas={xshift = 0.4cm}]
		    \draw[ultra thick,black,->] (S-1-8.north east) --(S-8-8.south east); 
  \end{scope}
      
\end{tikzpicture}
\caption{How to apply the SWM to the time series made up of the sequence of the 8 numerical values specified; see indications below.}
\label{f2}
\end{figure}

Consider the case for which $n=8$, specified in figure~\ref{f1}.  Figure~\ref{f2} illustrates that case.

Note that figure~\ref{f2} is composed of 8 rows and 8 columns. It may be considered to be an 8 by 8 matrix.

Look first at the eighth row, corresponding to the eighth train of square waves $S_8$ to be analyzed. In that row, the sequence $\dashline C_8\dashline {-C}_{8} \dashline$ corresponds to one square wave in $S_8$; $C_8$ corresponds to the first semi-wave of that square wave, and $-C_8$ to the second square wave. Which one of those semi-waves is positive and which is negative depends on the value to be computed (as specified below) for $C_8$. If $C_8$ is positive, then the first semi-wave will be positive, and the second will be negative. If, however, $C_8$ is negative, then the first semi-wave will be negative, and the second will be positive.

Note that in the last row, that of $S_8$, there are exactly 4 square waves; that is, those 4 square waves \textit{fit} in \SI{2}{\second} (the total duration of the interval $\Delta t$ during which time the values of the time series analyzed were measured). The frequency $f_8$ corresponding to that train of square waves $S_8$ is equal to the number of square waves that  \textit{fit} in the time unit \SI{1}{\second}; therefore, $f_8=\SI{2}{\hertz}$.

The wavelength of the next to the last train of square waves $S_7$ (corresponding to the next to the last line in figure~\ref{f2}), is twice the wavelength of the train of square waves $S_8$. The wavelength of the second to the last train of square waves $S_6$ is triple the wavelength of $S_8$, and so on successively. Hence, the wavelength $S_1$, the train of square waves in the first row in figure~\ref{f2}, is 8 times greater than the wavelength of $S_8$.

Given that the wavelength in $S_7$ is twice that of $S_8$, the frequency $f_7$ corresponding to $S_7$ is half that of $f_8$: $f_7=\tfrac{1}{2}f_8$. In figure~\ref{f2}, it can be seen that one square wave in $S_7$ \textit{fits} exactly in the unit of time \SI{1}{\second}; one square wave in $S_7$ is represented in figure~\ref{f2} as: $\dashline\;\; C_7\dashline\;\; C_7 \dashline -C_7 \dashline -C_7 \dashline $. 

To understand the formalism used, one may carefully verify at least some of the above results. Thus, for example, consider the following equation: $f_5=\tfrac{1}{4}f_8 = \tfrac{1}{4} \cdot \SI{2}{\hertz}= \SI[fraction-function=\tfrac,quotient-mode = fraction]{1/2}{\hertz} $. For the train of square waves $S_5$, one square wave can be represented in figure~\ref{f2} as:
$\scriptsize \dashline \; C_5 \dashline \quad C_5 \dashline \quad C_5 \dashline \quad C_5 \dashline {-C}_{5}\dashline {-C}_{5} \dashline {-C}_{5} \dashline {-C}_{5} \dashline $. For $S_5$, exactly one square wave \textit{fits} in the entire lapse in $\Delta t = \SI{2}{\second} $, as is clearly shown in figure~\ref{f2}; and exactly one semi-wave in $S_5$ \textit{fits} in the unit of time \SI{1}{\second}. So by definition of the notion of frequency, the number of waves per unit of time, the frequency $f_5$ corresponding to $S_5$ is: $f_5=\tfrac{\tfrac{1}{2}}{\SI{1}{\second}}=\tfrac{1}{\SI{2}{\second}} = \SI[fraction-function=\tfrac,quotient-mode = fraction]{1/2}{\hertz} $.

In the example discussed, 4 values were measured per second. In other words, for the sampling frequency $f_s$, this equation is valid: $f_s=\SI{4}{\hertz}$; the following relation between $f_8$ and $f_s$ is also valid: $f_8=\tfrac{1}{2}f_s$.

The results obtained for the frequencies $f_i$, where $i=1,2,\ldots, 8$, corresponding to the different trains of square waves $S_i$ are:
\begin{alignat*}{2}
f_1 &= \SI[fraction-function=\frac,quotient-mode = fraction]{2/8}{\hertz} &= \SI[fraction-function=\frac,quotient-mode = fraction]{1/4}{\hertz} \\
f_2 &= \SI[fraction-function=\frac,quotient-mode = fraction]{2/7}{\hertz} & \\
f_3 &= \SI[fraction-function=\frac,quotient-mode = fraction]{2/6}{\hertz} &= \SI[fraction-function=\frac,quotient-mode = fraction]{1/3}{\hertz} \\
f_4 &= \SI[fraction-function=\frac,quotient-mode = fraction]{2/5}{\hertz} &  \\
f_5 &= \SI[fraction-function=\frac,quotient-mode = fraction]{2/4}{\hertz} &= \SI[fraction-function=\frac,quotient-mode = fraction]{1/2}{\hertz} \\
f_6 &= \SI[fraction-function=\frac,quotient-mode = fraction]{2/3}{\hertz} &  \\
f_7 &= \SI[fraction-function=\frac,quotient-mode = fraction]{2/2}{\hertz} &= \SI{1}{\hertz} \\
f_8 &= \SI[fraction-function=\frac,quotient-mode = fraction]{2/1}{\hertz} &= \SI{2}{\hertz}
\end{alignat*}

How to compute the values of the 8 coefficients $C_1$, $C_2$, $C_3$,  $C_4$, $C_5$, $C_6$, $C_7$, and $C_8$, displayed in figure~\ref{f2} will be specified below.

In this example, the interval $\Delta t=\SI{2}{\second}$ was divided into 8 equal subintervals because the time series analyzed is made up of exactly 8 equal values. If that time series were composed of $50,000$ values, then the interval $\Delta t$ would be divided into 50,000 equal subintervals. In general, if the time series to be analyzed is made up of $n$ values, where $n=1,2,3,\ldots$, then $\Delta t$ is divided into $n$ equal subintervals.

First, look at the left column in figure~\ref{f2}, and add up the elements in that column, following the order indicated by the vertical arrow pointing down (at the right of the figure) and the result is made equal to the first value, \SI{84}{\mV}, of that time series. Hence, the following equation is obtained:
\begin{equation*}
C_1+C_2+C_3+C_4+C_5+C_6+C_7+C_8=\SI{84}{\mV}
\end{equation*}

Note that $C_1$, $C_2$, $C_3$,  $C_4$, $C_5$, $C_6$, $C_7$, and $C_8$ are respectively the values of the trains of square waves $S_1$, $S_2$, $S_3$,  $S_4$, $S_5$, $S_6$, $S_7$, and $S_8$, at the midpoint of the first subinterval of $\Delta t$.

Then, look at the second column in figure~\ref{f2}, and add up the elements in that column making the result equal to the second value, \SI{-152}{\mV}, in that time series. Thus, the following equation is obtained:
\begin{equation*}
C_1+C_2+C_3+C_4+C_5+C_6+C_7-C_8=\SI{-152}{\mV}
\end{equation*}

Observe that $C_1$, $C_2$, $C_3$,  $C_4$, $C_5$, $C_6$, $C_7$, and $-C_8$ are respectively the values of the trains of square waves $S_1$, $S_2$, $S_3$,  $S_4$, $S_5$, $S_6$, $S_7$, and $S_8$, at the midpoint of the second subinterval of $\Delta t$.

Now look at the third column in figure~\ref{f2}, and add up the elements in that column, making the result equal to the third value, \SI{63}{\mV}, in that time series. The following equation is obtained:
\begin{equation*}
C_1+C_2+C_3+C_4+C_5+C_6-C_7+C_8=\SI{63}{\mV}
\end{equation*}

Note that $C_1$, $C_2$, $C_3$,  $C_4$, $C_5$, $C_6$, $-C_7$, and $C_8$ are respectively the values of the trains of square waves $S_1$, $S_2$, $S_3$,  $S_4$, $S_5$, $S_6$, $S_7$, and $S_8$, at the midpoint of the third subinterval of $\Delta t$.

The fourth, fifth, sixth and seventh equations are obtained using the same type of procedure.

Finally, that procedure is applied to the eighth column in ~\ref{f2}, and the following equation is obtained:
\begin{equation*}
C_1-C_2-C_3-C_4-C_5+C_6-C_7-C_8=\SI{-14}{\mV}
\end{equation*}

Observe that $C_1$, $-C_2$, $-C_3$,  $-C_4$, $-C_5$, $C_6$, $-C_7$, and $-C_8$ are respectively the values of the trains of square waves for $S_1$, $S_2$, $S_3$,  $S_4$, $S_5$, $S_6$, $S_7$, and $S_8$, at the midpoint of the eighth subinterval of $\Delta t$.

The following system of 8 linear algebraic equations was obtained as specified above.
\begin{equation} 
  \left.\begin{aligned}
C_1+C_2+C_3+C_4+C_5+C_6+C_7+C_8&=&\SI{84}{\mV} \\
C_1+C_2+C_3+C_4+C_5+C_6+C_7-C_8&=&\SI{-152}{\mV}\\
C_1+C_2+C_3+C_4+C_5+C_6-C_7+C_8&=&\SI{63}{\mV}\\
C_1+C_2+C_3+C_4+C_5-C_6-C_7-C_8&=&\SI{98}{\mV}\\
C_1+C_2+C_3+C_4-C_5-C_6+C_7+C_8&=&\SI{-35}{\mV}\\
C_1+C_2+C_3-C_4-C_5-C_6+C_7-C_8&=&\SI{0}{\mV}\\
C_1+C_2-C_3-C_4-C_5+C_6-C_7+C_8&=&\SI{145}{\mV}\\
C_1-C_2-C_3-C_4-C_5+C_6-C_7-C_8&=&\SI{-14}{\mV}
  \end{aligned}
  \qquad \right\}
  \label{eq1}
\end{equation}

The above system of 8 linear algebraic equations~(\ref{eq1}) has 8 unknowns: $C_1$, $C_2$, $C_3$,  $C_4$, $C_5$, $C_6$, $C_7$, and $C_8$, which are precisely the 8 coefficients whose values must be computed.
 
Once that system of equations~(\ref{eq1}) has been solved, the values for the coefficients are: 
\begin{equation*}
\begin{aligned}
 C_1 &= &\SI{170.5}{\mV} & \quad  & C_5 &= &\SI{195.0}{\mV} \\
 C_2 &= &\SI{-38.5}{\mV}&   & C_6 &=&\SI{-135.5}{\mV}\\
 C_3 &= &\SI{-100.5}{\mV} &  & C_7 &= &\SI{10.5}{\mV}\\
 C_4 &=& \SI{-135.5}{\mV}&   & C_8 &=& \SI{118.0}{\mV} \\
\end{aligned} 
\end{equation*}

If in each member on the left side of the system of equations (\ref{eq1}), the coefficients $C_1$, $C_2$, $C_3$,  $C_4$, $C_5$, $C_6$, $C_7$, and $C_8$ are replaced by the values computed for them, and if the resulting 8 algebraic sums are completed, the following system of 8 equations is obtained:
\begingroup\makeatletter\def\f@size{9.5}\check@mathfonts
 \begin{equation}
  \left.\begin{aligned}
(170.5-38.5-100.5-135.5+195-135.5+10.5+118)\,\si{\mV}&=&\SI{84}{\mV} \\
(170.5-38.5-100.5-135.5+195-135.5+10.5-118)\,\si{\mV}&=&\SI{-152}{\mV}\\
(170.5-38.5-100.5-135.5+195-135.5-10.5+118)\,\si{\mV}&=&\SI{63}{\mV}\\
(170.5-38.5-100.5-135.5+195+135.5-10.5-118)\,\si{\mV}&=&\SI{98}{\mV}\\
(170.5-38.5-100.5-135.5-195+135.5+10.5+118)\,\si{\mV}&=&\SI{-35}{\mV}\\
(170.5-38.5-100.5+135.5-195+135.5+10.5-118)\,\si{\mV}&=&\SI{0}{\mV}\\
(170.5-38.5+100.5+135.5-195-135.5-10.5+118)\,\si{\mV}&=&\SI{145}{\mV}\\
(170.5+38.5+100.5+135.5-195-135.5-10.5-118)\,\si{\mV}&=&\SI{-14}{\mV}
  \end{aligned}
   \right\}
  \label{eq2}
\end{equation}
\endgroup
Note that the members on the right of (\ref{eq2}) coincide exactly with the members on the right of (\ref{eq1}).

In figure~\ref{f3} the parts of the trains of square waves  $S_1$, $S_2$, $S_3$,  $S_4$, $S_5$, $S_6$, $S_7$, $S_8$ corresponding to interval $\Delta t$ have been displayed.

\setcounter{figure}{2}
\begin{figure}[H]
\centering
\subfloat[$S_1(t)$.]{
\begin{tikzpicture} [x=40mm,y=.1mm]
\draw[->](0, 0)--(2.2, 0) ;
\draw[->](0, -250)--(0, 250);
\foreach \x in { 1,2}
     		\draw(\x, 1pt)--(\x, -3pt)
            node[below right=-2pt] {\x};
\foreach \y in { -200,-100,...,200}
     		\draw(1pt,\y)--(-3pt,\y)
                 node[anchor = east] {\y};
\node at (-.15, 310) {\si{\mV}};
\node at (-0.15,260) {$S_1(t)$}; 
\node at(2.22,-25) {$t$(\si{\second})};  
\draw[thick](0,170.5)--(2,170.5);
\draw[dotted](2,170.5)--(2,-170.5);
\fill (0.125,170.5) circle [radius=2pt];
\fill (0.375,170.5) circle [radius=2pt];
\fill (0.625,170.5) circle [radius=2pt];
\fill (0.875,170.5) circle [radius=2pt];
\fill (1.125,170.5) circle [radius=2pt];
\fill (1.375,170.5) circle [radius=2pt];
\fill (1.625,170.5) circle [radius=2pt];
\fill (1.875,170.5) circle [radius=2pt];
\draw[dashed] (0.25,-200) -- (0.25,200);
\draw[dashed] (0.5,-200) -- (0.5,200);
\draw[dashed] (0.75,-200) -- (0.75,200);
\draw[dashed] (1.0,-200) -- (1.0,200);
\draw[dashed] (1.25,-200) -- (1.25,200);
\draw[dashed] (1.50,-200) -- (1.5,200);
\draw[dashed] (1.75,-200) -- (1.75,200);
\draw[dashed] (2,-200) -- (2,200);
\end{tikzpicture}
\label{f3a}
}\qquad
\phantomcaption
\end{figure}
\begin{figure}
\ContinuedFloat
\subfloat[$S_2(t)$.]{
\begin{tikzpicture} [x=40mm,y=.1mm]
\draw[->](0, 0)--(2.2, 0) ;
\draw[->](0, -250)--(0, 250);
\foreach \x in { 1,2}
     		\draw(\x, 1pt)--(\x, -3pt)
            node[below right=-2pt] {\x};
\foreach \y in { -200,-100,...,200}
     		\draw(1pt,\y)--(-3pt,\y)
                 node[anchor = east] {\y};
\node at (-.15, 310) {\si{\mV}};
\node at (-0.15,260) {$S_1(t)$}; 
\node at(2.22,-25) {$t$(\si{\second})}; 
\draw[thick](0,-38.5)--(1.7500000000000000000000000001,-38.5);
\draw[thick](1.7500000000000000000000000001,38.5)--(2.0000000000000000000000000002,38.5);
\draw[dotted](1.7500000000000000000000000001,-38.5)--(1.7500000000000000000000000001,38.5);
\fill (0.125,-38.5) circle [radius=2pt];
\fill (0.375,-38.5) circle [radius=2pt];
\fill (0.625,-38.5) circle [radius=2pt];
\fill (0.875,-38.5) circle [radius=2pt];
\fill (1.125,-38.5) circle [radius=2pt];
\fill (1.375,-38.5) circle [radius=2pt];
\fill (1.625,-38.5) circle [radius=2pt];
\fill (1.875,38.5) circle [radius=2pt];
\draw[dashed] (0.25,-200) -- (0.25,200);
\draw[dashed] (0.5,-200) -- (0.5,200);
\draw[dashed] (0.75,-200) -- (0.75,200);
\draw[dashed] (1.0,-200) -- (1.0,200);
\draw[dashed] (1.25,-200) -- (1.25,200);
\draw[dashed] (1.50,-200) -- (1.5,200);
\draw[dashed] (1.75,-200) -- (1.75,200);
\draw[dashed] (2,-200) -- (2,200);
\end{tikzpicture}
\label{f3b}
}\qquad
\subfloat[$S_3(t)$.]{
\begin{tikzpicture} [x=40mm,y=.1mm]
\draw[->](0, 0)--(2.2, 0) ;
\draw[->](0, -250)--(0, 250);
\foreach \x in { 1,2}
     		\draw(\x, 1pt)--(\x, -3pt)
            node[below right=-2pt] {\x};
\foreach \y in { -200,-100,...,200}
     		\draw(1pt,\y)--(-3pt,\y)
                 node[anchor = east] {\y};
\node at (-.15, 310) {\si{\mV}};
\node at (-0.15,260) {$S_3(t)$}; 
\node at(2.22,-25) {$t$(\si{\second})}; 
\draw[thick](0,-100.5)--(1.5000000000000000000000000002,-100.5);
\draw[thick](1.5000000000000000000000000002,100.5)--(2.0000000000000000000000000004,100.5);
\draw[dotted](1.5000000000000000000000000002,-100.5)--(1.5000000000000000000000000002,100.5);
\fill (0.125,-100.5) circle [radius=2pt];
\fill (0.375,-100.5) circle [radius=2pt];
\fill (0.625,-100.5) circle [radius=2pt];
\fill (0.875,-100.5) circle [radius=2pt];
\fill (1.125,-100.5) circle [radius=2pt];
\fill (1.375,-100.5) circle [radius=2pt];
\fill (1.625,100.5) circle [radius=2pt];
\draw[dashed] (0.25,-200) -- (0.25,200);
\draw[dashed] (0.5,-200) -- (0.5,200);
\draw[dashed] (0.75,-200) -- (0.75,200);
\draw[dashed] (1.0,-200) -- (1.0,200);
\draw[dashed] (1.25,-200) -- (1.25,200);
\draw[dashed] (1.50,-200) -- (1.5,200);
\draw[dashed] (1.75,-200) -- (1.75,200);
\draw[dashed] (2,-200) -- (2,200);
\fill (1.875,100.5) circle [radius=2pt];
\end{tikzpicture}
\label{f3c}
}\qquad
\phantomcaption
\end{figure}
\begin{figure}
\ContinuedFloat
\subfloat[$S_4(t)$.]{
\begin{tikzpicture} [x=40mm,y=.1mm]
\draw[->](0, 0)--(2.2, 0) ;
\draw[->](0, -250)--(0, 250);
\foreach \x in { 1,2}
     		\draw(\x, 1pt)--(\x, -3pt)
            node[below right=-2pt] {\x};
\foreach \y in { -200,-100,...,200}
     		\draw(1pt,\y)--(-3pt,\y)
                 node[anchor = east] {\y};
\node at (-.15, 310) {\si{\mV}};
\node at (-0.15,260) {$S_4(t)$}; 
\node at(2.22,-25) {$t$(\si{\second})};
\draw[thick](0,-135.5)--(1.25,-135.5);
\draw[thick](1.25,135.5)--(2.00,135.5);
\draw[dotted](1.25,-135.5)--(1.25,135.5);
\fill (0.125,-135.5) circle [radius=2pt];
\fill (0.375,-135.5) circle [radius=2pt];
\fill (0.625,-135.5) circle [radius=2pt];
\fill (0.875,-135.5) circle [radius=2pt];
\fill (1.125,-135.5) circle [radius=2pt];
\fill (1.375,135.5) circle [radius=2pt];
\fill (1.625,135.5) circle [radius=2pt];
\fill (1.875,135.5) circle [radius=2pt];
\draw[dashed] (0.25,-200) -- (0.25,200);
\draw[dashed] (0.5,-200) -- (0.5,200);
\draw[dashed] (0.75,-200) -- (0.75,200);
\draw[dashed] (1.0,-200) -- (1.0,200);
\draw[dashed] (1.25,-200) -- (1.25,200);
\draw[dashed] (1.50,-200) -- (1.5,200);
\draw[dashed] (1.75,-200) -- (1.75,200);
\draw[dashed] (2,-200) -- (2,200);
\end{tikzpicture}
\label{f3d}
}\qquad
\subfloat[$S_5(t)$.]{
\begin{tikzpicture} [x=40mm,y=.1mm]
\draw[->](0, 0)--(2.2, 0) ;
\draw[->](0, -250)--(0, 250);
\foreach \x in { 1,2}
     		\draw(\x, 1pt)--(\x, -3pt)
            node[below right=-2pt] {\x};
\foreach \y in { -200,-100,...,200}
     		\draw(1pt,\y)--(-3pt,\y)
                 node[anchor = east] {\y};
\node at (-.15, 310) {\si{\mV}};
\node at (-0.15,260) {$S_5(t)$}; 
\node at(2.22,-25) {$t$(\si{\second})};
\draw[thick](0,195.0)--(1,195.0);
\draw[thick](1,-195.0)--(2,-195.0);
\draw[dotted](1,195.0)--(1,-195.0);
\draw[dotted](2,195.0)--(2,-195.0);
\fill (0.125,195.0) circle [radius=2pt];
\fill (0.375,195.0) circle [radius=2pt];
\fill (0.625,195.0) circle [radius=2pt];
\fill (0.875,195.0) circle [radius=2pt];
\fill (1.125,-195.0) circle [radius=2pt];
\fill (1.375,-195.0) circle [radius=2pt];
\fill (1.625,-195.0) circle [radius=2pt];
\fill (1.875,-195.0) circle [radius=2pt];
\draw[dashed] (0.25,-200) -- (0.25,200);
\draw[dashed] (0.5,-200) -- (0.5,200);
\draw[dashed] (0.75,-200) -- (0.75,200);
\draw[dashed] (1.0,-200) -- (1.0,200);
\draw[dashed] (1.25,-200) -- (1.25,200);
\draw[dashed] (1.50,-200) -- (1.5,200);
\draw[dashed] (1.75,-200) -- (1.75,200);
\draw[dashed] (2,-200) -- (2,200);
\end{tikzpicture}
\label{f3e}
}\qquad
\phantomcaption
\end{figure}
\begin{figure}
\ContinuedFloat
\subfloat[$S_6(t)$.]{
\begin{tikzpicture} [x=40mm,y=.1mm]
\draw[->](0, 0)--(2.2, 0) ;
\draw[->](0, -250)--(0, 250);
\foreach \x in { 1,2}
     		\draw(\x, 1pt)--(\x, -3pt)
            node[below right=-2pt] {\x};
\foreach \y in { -200,-100,...,200}
     		\draw(1pt,\y)--(-3pt,\y)
                 node[anchor = east] {\y};
\node at (-.15, 310) {\si{\mV}};
\node at (-0.15,260) {$S_6(t)$}; 
\node at(2.22,-25) {$t$(\si{\second})};
\draw[thick](0,-135.5)--(0.75,-135.5);
\draw[thick](0.75,135.5)--(1.50,135.5);
\draw[dotted](0.75,-135.5)--(0.75,135.5);
\draw[dotted](1.50,-135.5)--(1.50,135.5);
\fill (0.125,-135.5) circle [radius=2pt];
\fill (0.375,-135.5) circle [radius=2pt];
\fill (0.625,-135.5) circle [radius=2pt];
\fill (0.875,135.5) circle [radius=2pt];
\fill (1.125,135.5) circle [radius=2pt];
\fill (1.375,135.5) circle [radius=2pt];
\draw[thick](1.50,-135.5)--(2.00,-135.5);
\draw[dotted](2.25,-135.5)--(2.25,135.5);
\draw[dotted](3.00,-135.5)--(3.00,135.5);
\fill (1.625,-135.5) circle [radius=2pt];
\fill (1.875,-135.5) circle [radius=2pt];
\draw[dashed] (0.25,-200) -- (0.25,200);
\draw[dashed] (0.5,-200) -- (0.5,200);
\draw[dashed] (0.75,-200) -- (0.75,200);
\draw[dashed] (1.0,-200) -- (1.0,200);
\draw[dashed] (1.25,-200) -- (1.25,200);
\draw[dashed] (1.50,-200) -- (1.5,200);
\draw[dashed] (1.75,-200) -- (1.75,200);
\draw[dashed] (2,-200) -- (2,200);
\end{tikzpicture}
\label{f3f}
}\qquad
\subfloat[$S_7(t)$.]{
\begin{tikzpicture} [x=40mm,y=.1mm]
\draw[->](0, 0)--(2.2, 0) ;
\draw[->](0, -250)--(0, 250);
\foreach \x in { 1,2}
     		\draw(\x, 1pt)--(\x, -3pt)
            node[below right=-2pt] {\x};
\foreach \y in { -200,-100,...,200}
     		\draw(1pt,\y)--(-3pt,\y)
                 node[anchor = east] {\y};
\node at (-.15, 310) {\si{\mV}};
\node at (-0.15,260) {$S_7(t)$}; 
\node at(2.22,-25) {$t$(\si{\second})};
\draw[thick](0,10.5)--(0.5,10.5);
\draw[thick](0.5,-10.5)--(1.0,-10.5);
\draw[dotted](0.5,10.5)--(0.5,-10.5);
\draw[dotted](1.0,10.5)--(1.0,-10.5);
\fill (0.125,10.5) circle [radius=2pt];
\fill (0.375,10.5) circle [radius=2pt];
\fill (0.625,-10.5) circle [radius=2pt];
\fill (0.875,-10.5) circle [radius=2pt];
\draw[thick](1.0,10.5)--(1.5,10.5);
\draw[thick](1.5,-10.5)--(2.0,-10.5);
\draw[dotted](1.5,10.5)--(1.5,-10.5);
\draw[dotted](2.0,10.5)--(2.0,-10.5);
\fill (1.125,10.5) circle [radius=2pt];
\fill (1.375,10.5) circle [radius=2pt];
\fill (1.625,-10.5) circle [radius=2pt];
\fill (1.875,-10.5) circle [radius=2pt];
\draw[dashed] (0.25,-200) -- (0.25,200);
\draw[dashed] (0.5,-200) -- (0.5,200);
\draw[dashed] (0.75,-200) -- (0.75,200);
\draw[dashed] (1.0,-200) -- (1.0,200);
\draw[dashed] (1.25,-200) -- (1.25,200);
\draw[dashed] (1.50,-200) -- (1.5,200);
\draw[dashed] (1.75,-200) -- (1.75,200);
\draw[dashed] (2,-200) -- (2,200);
\end{tikzpicture}
\label{f3g}
}\qquad
\phantomcaption
\end{figure}
\begin{figure}
\ContinuedFloat
\subfloat[$S_8(t)$.]{
\begin{tikzpicture} [x=40mm,y=.1mm]
\draw[->](0, 0)--(2.2, 0) ;
\draw[->](0, -250)--(0, 250);
\foreach \x in { 1,2}
     		\draw(\x, 1pt)--(\x, -3pt)
            node[below right=-2pt] {\x};
\foreach \y in { -200,-100,...,200}
     		\draw(1pt,\y)--(-3pt,\y)
                 node[anchor = east] {\y};
\node at (-.15, 310) {\si{\mV}};
\node at (-0.15,260) {$S_8(t)$}; 
\node at(2.22,-25) {$t$(\si{\second})}; 
\draw[thick](0,118.0)--(0.25,118.0);
\draw[thick](0.25,-118.0)--(0.50,-118.0);
\draw[dotted](0.25,118.0)--(0.25,-118.0);
\draw[dotted](0.50,118.0)--(0.50,-118.0);
\fill (0.125,118.0) circle [radius=2pt];
\fill (0.375,-118.0) circle [radius=2pt];
\draw[thick](0.50,118.0)--(0.75,118.0);
\draw[thick](0.75,-118.0)--(1.00,-118.0);
\draw[dotted](0.75,118.0)--(0.75,-118.0);
\draw[dotted](1.00,118.0)--(1.00,-118.0);
\fill (0.625,118.0) circle [radius=2pt];
\fill (0.875,-118.0) circle [radius=2pt];
\draw[thick](1.00,118.0)--(1.25,118.0);
\draw[thick](1.25,-118.0)--(1.50,-118.0);
\fill (1.125,118.0) circle [radius=2pt];
\fill (1.375,-118.0) circle [radius=2pt];
\draw[thick](1.50,118.0)--(1.75,118.0);
\draw[thick](1.75,-118.0)--(2.00,-118.0);
\draw[dotted](1.75,118.0)--(1.75,-118.0);
\draw[dotted](2.00,118.0)--(2.00,-118.0);
\fill (1.625,118.0) circle [radius=2pt];
\fill (1.875,-118.0) circle [radius=2pt];
\draw[dashed] (0.25,-200) -- (0.25,200);
\draw[dashed] (0.5,-200) -- (0.5,200);
\draw[dashed] (0.75,-200) -- (0.75,200);
\draw[dashed] (1.0,-200) -- (1.0,200);
\draw[dashed] (1.25,-200) -- (1.25,200);
\draw[dashed] (1.50,-200) -- (1.5,200);
\draw[dashed] (1.75,-200) -- (1.75,200);
\draw[dashed] (2,-200) -- (2,200);
\end{tikzpicture}
\label{f3h}
}
\caption{Trains of square waves $S_1, S_2, S_3, \dots$, and $S_8$}
\label{f3}
\end{figure}
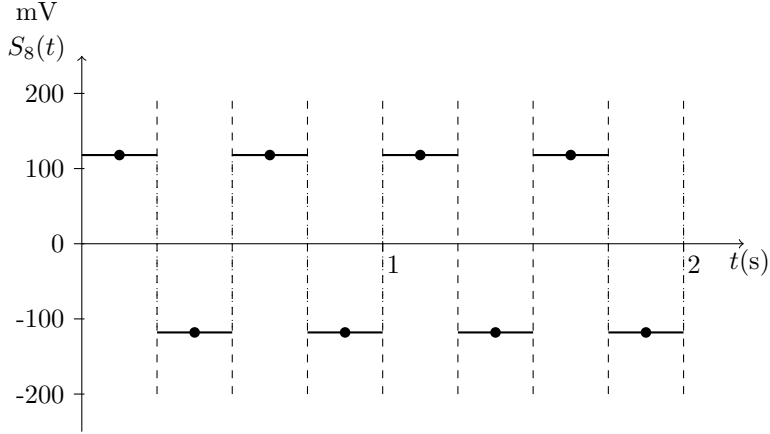
\newpage

\setcounter{figure}{3}

For each of the 8 equal subintervals of $\Delta t$, the algebraic sum of the values highlighted in those trains of square waves must be computed. The algebraic sum corresponding to the $i$-th subinterval ($i=1,2,3,\ldots, 8$) must be equal to the $i$-th value of the time series considered. 

In general, this approach can be applied to analyze any time series composed of $n$ measured values, such that $n=1,2,3,\ldots$, during a specific time interval $\Delta t$ expressed in seconds. It is accepted that the lapse between two consecutive values in the time series analyzed will remain constant. 

The frequency of the measurement of values in the time series (or the \textit{sampling frequency} $f_s$) is the number of values measured per second (the time unit). Therefore, the number $n$ of measured values in a time interval $\Delta t$ can be expressed as:
\begin{equation*}
n=f_s\cdot \Delta t.
 \end{equation*}
 
If the time series is composed of $n$ values, the interval $\Delta t$ must be divided into $n$ subintervals of equal duration. 

To analyze that time series with the SWM, it is essential to construct an $n \times n$-matrix like that of figure~\ref{f2}.

The last row in that matrix corresponds to the \textit{n}-th train of square waves $S_n$ of the $n$ trains of square waves to be considered in this particular case.

Each square wave in $S_n$ can be symbolized as $\dashline  {C}_{n}\dashline {-C}_{n} \dashline$, occupying 2 subintervals of the $n$ equal subintervals into which $\Delta t$ was divided. Thus, the number of square waves in $S_n$ per time unit ($f_n$) can be expressed as follows:

\begin{equation*}
f_n=\frac{1}{2} \cdot \frac{n}{\Delta t}
 \end{equation*}

The next to the last row of this $n \times n$-matrix is that of the train of square waves $S_{n-1}$; each square wave is represented as $\dashline {C}_{n-1} \dashline {C}_{n-1} \allowbreak \dashline {-C}_{n-1} \allowbreak \dashline {-C}_{n-1} \dashline$. It can be seen that the length of each of the square waves composing $S_{n-1}$ is twice that of the wave of each square wave composing $S _{n}$. Hence, the frequency $f_{n-1}$ corresponding to $S_{n-1}$ is equal to half of frequency $f_{n}$ corresponding to $S_{n}$:
\begin{equation*}
f_{n-1}=\frac{1}{2} \cdot f_n
 \end{equation*}

The second to the last row in this $n \times n$-matrix corresponds to the train of square waves $S_{n-2}$. Each square wave in $S_{n-2}$ can be represented as $\dashline {C}_{n-2} \allowbreak  \dashline {C}_{n-2} \dashline {C}_{n-2} \dashline {-C}_{n-2}\dashline {-C}_{n-2}\dashline {-C}_{n-2} \dashline$. It can be observed that the length of each of the square waves composing $S_{n-2}$ is three times that of the wave of each square wave composing $S _{n}$. Thus, the frequency $f_{n-2}$ corresponding to $S_{n-2}$ is equal to a third of the frequency $f_{n}$ corresponding to $S_{n}$:  
\begin{equation*}
f_{n-2}=\frac{1}{3} \cdot f_n
 \end{equation*}

In general, the frequency $f_{i}$ corresponding to the $i$-th train of square waves $S_{i}$ will be expressed as follows:
\begin{equation}
f_{i}=f_{n-(n-i)} = \left(\frac{1}{n-i+1}\right)\cdot f_n =\left(\frac{1}{n-i+1}\right)\cdot\frac{1}{2}\cdot\frac{n}{\Delta t} = \frac{1}{2\Delta t}\cdot\frac{n}{n-i+1}
\label{eq3}
 \end{equation}
 
The values of the coefficients $C_1,C_2,\ldots, C_n$ are computed with a procedure like that used to compute $C_1,C_2,\ldots, C_8$, in the example analyzed of a time series composed of a sequence of 8 values. 
 
To begin with, the algebraic sum of all the elements making up the first column of the $n \times n$-matrix is made equal to the first value of the time series to be analyzed. 

Next, the algebraic sum of all the elements composing the second column of the $n \times n$-matrix is made equal to the second value of the time series to be analyzed. 

And so on up to the algebraic sum of all the elements of the $n$-th column of the $n \times n$-matrix, which is made equal to the last value ($n$-th value) of the time series to be analyzed. 

In this way, a system of $n$ linear algebraic equations is obtained. The coefficients $C_1,C_2,\ldots, C_n$, which are the unknowns of that system of equations, can be computed. 
 
An algorithm to make it easier to determine the system of linear algebraic equations will be described in section 3. This system of equations must be solved whenever the SWM is used to analyze a time series.
 
 Suppose that the SWM is being used to analyze a time series made up of a sequence of $n$ measured values of a signal depending on time. As explained above, each train of square waves $S_i$, for $i=1,2,\ldots,n$, corresponding to that analysis, will be characterized by the value of a frequency $f_i$ and the value of a coefficient $C_i$ that can be computed using the procedures described above. Thus the result of the analysis can be expressed by a sequence of $n$ dyads. The first element of the first of those dyads is the value of $f_1$, and the second element of the first dyad is $C_1$. The first element of the second dyad is the value of $f_2$, and the second element of that second dyad is $C_2$, and so on, successively; hence, the first element of the $n$-th dyad is the value of $f_n$, and the second element of that dyad is $C_n$.
 
Thus, for example, the 8 dyads that make  it possible to express the result of the analysis carried out with the SWM of the sequence of 8 values represented in figure~\ref{f1} are:
\begin{equation*}
\begin{aligned}[c]
(f_{1}; C_{1}) &=(0.250000; 170.500000)  \\
(f_{3}; C_{3}) &=(0.333333; -100.500000)   \\
(f_{5}; C_{5}) &= (0.500000; 195.000000)  \\
(f_{7}; C_{7}) &= (1.000000; 10.500000)\\
\end{aligned}
\qquad
\begin{aligned}[c]
(f_{2}; C_{2})&=   (0.285714; -38.500000)\\
(f_{4}; C_{4}) &=    (0.400000; -135.500000)  \\
(f_{6}; C_{6}) &=   (0.666667; -135.500000)  \\
(f_{8}; C_{8}) &=   (2.000000; 118.000000) \\
\end{aligned}
\end{equation*}

A mathematical tool (the Square Wave Transform, SWT) was introduced in \cite{b0} to make it possible to represent the frequency domain of the results of any analysis carried out with the SWM. The values of $f_i$, for $i=1,2,\ldots, n$, are plotted on the x-axis of an orthogonal Cartesian coordinate system. For every $f_i$, the value of the corresponding $C_i$ is represented by a vertical bar whose value is seen on the y-axis.

\begin{figure}[H]
\centering
\includegraphics[width=4.5in]{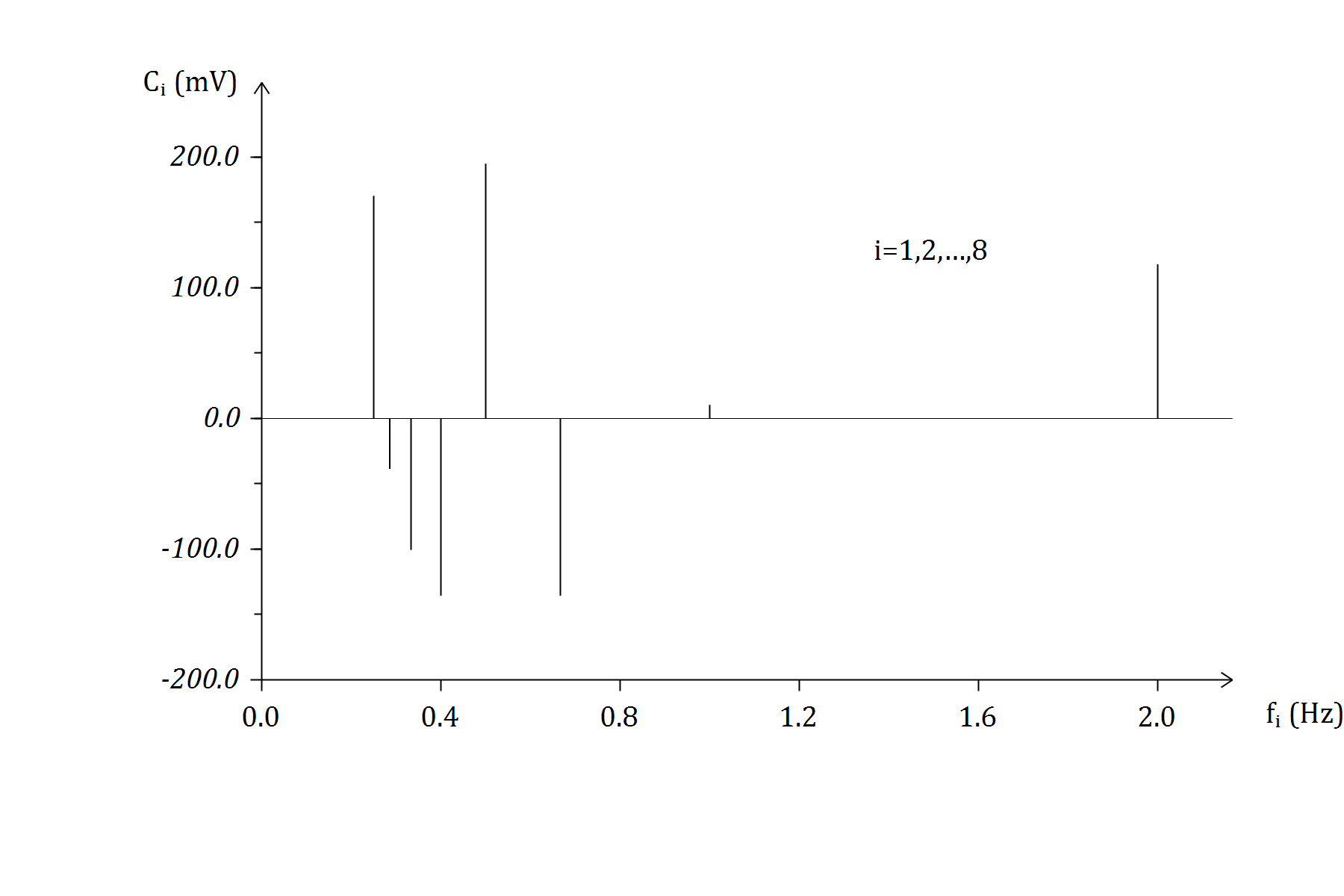}
\caption{SWT corresponding to the time series composed of the sequence of 8 voltage measurements given in section 2}
\label{f4}
\end{figure}

The SWT corresponding to the sequence of 8 voltage values specified in figure \ref{f1} has been represented in figure \ref{f4} to display the information that was provided quantitatively in more detail by the sequence of 8 dyads specified above.
 
\section{Description of an algorithm to facilitate the determination of the required system of linear algebraic equations}

Based on the regularities present in the tables of coefficients such as in figure 2, the algorithm described here may be used to determine the system of linear algebraic equations to be solved whenever the SWM is used to analyze a time series.

For this purpose, let $Q$ be the integer quotient in a division operation, and $R$ the integer remainder, as exemplified: For $10 \div 7$, $Q=1$ and $R=3$; or for $7 \div 10$, $Q=0$ and $R=7$.

In general, if $N_1$ and $N_2$, are two natural numbers, the quotient of $N_1$ and $N_2$ can be expressed in terms of $Q$ and $R$ as:

\begin{equation*}
\frac{N_1}{N_2}=Q+\frac{R}{N_2}.
\end{equation*}\\
If $N_1<N_2$, then $Q=0$ and $R=N_2$. These results will be used below. 

If one considers a time series of $n$ values (which can result from a sequence of $n$ measurements taken during a time interval $\Delta t$) such that between any two consecutive values of those $n$ values the duration is the same, then the system of $n$ linear algebraic equations that must be solved in order to analyze that time series with the SWM can be presented initially as follows: 

\begin{equation}
  \left.\begin{aligned}
      \sum_{j=1}^n  C_{1,j}&= V_1 \\
      \sum_{j=1}^n  C_{2,j}&= V_2 \\
            \sum_{j=1}^n  C_{3,j}&= V_3 \\
  \vdotswithin{=}\\
  \sum_{j=1}^n  C_{n,j}&= V_n \\
    \end{aligned}
  \qquad \right\}
  \label{eq4}
\end{equation}\\

$V_1, V_2,V_3,\ldots, V_n$ are the measured values corresponding to the first subinterval of the $n$ subintervals into which $\Delta t$ is divided, to the second subinterval, to the third subinterval,  $\ldots$, and the $n$-th subinterval of $\Delta t$ respectively.

Consider each of the coefficients $C_{i,j}$ in the above system of equations (\ref{eq4}). In the first place, attention is given to the second subscript ($j$) of that coefficient. Starting with  $j$, $l_j=n-j+1$ is computed. (Recall that $n$ is the number of subintervals into which $\Delta t$ was divided.) Thus, for example, if $n=8$ and $j=5$, then $l_j=8-5+1=4$.

Two other examples of how $l_j$ is computed are discussed below. If $n=10$ and $j=10$, then $l_j=10-10+1=1$. If $n=10$ and $j=1$, then $l_j=10-1+1=10$. 

Note that $l_j$ is the length of the semi-wave (that is, of the half-wave) corresponding to the train of square waves $S_j$, expressed in the number of subintervals in $\Delta t$.

For each coefficient $C_{i,j}$, consideration must be given to the values of $Q$ and $R$ resulting from the following division:
\begin{equation*}
\frac{i}{l_j}=\frac{i}{n-j+1}.
\end{equation*}

A variable $K$ is introduced for one of the following values: $0,1,2,3,\ldots$; thus the following result is obtained:
\begin{enumerate}
\item If $(Q=2K\ \mathrm{and}\  R=0)$, then $C_{i,j}=-C_j$
\item If $(Q=2K\ \mathrm{and}\  R\neq0)$, then $C_{i,j}=+C_j$
\item If $(Q=2K+1\ \mathrm{and}\  R=0)$, then $C_{i,j}=+C_j$
\item If $(Q=2K+1\ \mathrm{and}\  R\neq0)$, then $C_{i,j}=-C_j$
\end{enumerate}

As an example of the application of this algorithm, the system of eight linear algebraic equations from section 2 will be obtained once again. Initially, this system of equations can be presented as follows:

\begin{equation*}
\begin{aligned}
      \sum_{j=1}^8  C_{1,j}&=& \SI{84}{mV} \\
      \sum_{j=1}^8  C_{2,j}&=& \SI{-152}{mV} \\
      \sum_{j=1}^8  C_{3,j}&=&\SI{63}{mV} \\
      \sum_{j=1}^8  C_{4,j}&=& \SI{98}{mV} \\
      \sum_{j=1}^8  C_{5,j}&=& \SI{-35}{mV} \\
      \sum_{j=1}^8  C_{6,j}&= &\SI{0}{mV} \\
      \sum_{j=1}^8  C_{7,j}&= &\SI{145}{mV}  \\      
      \sum_{j=1}^8  C_{8,j}&=& \SI{-14}{mV} \\
    \end{aligned}
\end{equation*}
 
Note that in the above system of equations, the first subscript of each coefficient $i$, for $i=1,2,3,\ldots,8$, indicates the order of the equation itself in which that coefficient appears; that is, the first, second, third, $\ldots$ or eighth equation and the second coefficient $j$, for $j=1,2,3,\ldots,8$ specifies the order in which that coefficient is present in the given equation (first, second third, $\ldots$ or eighth). Thus, for example, $C_{5,4}$ is a coefficient present in the fifth of the above equations and it appears in the fourth place of that equation.

Each coefficient $C_{i,j}$ of the 64 coefficients in the above system of equations can then be replaced by $-C_j$ or by $+C_j$, as established in 1), 2), 3) and 4).

When proceeding in this way, the system of linear algebraic equations (1) mentioned in section 2 is obtained:

\begin{equation}
  \left.\begin{aligned}
C_1+C_2+C_3+C_4+C_5+C_6+C_7+C_8&=&\SI{84}{\mV} \\
C_1+C_2+C_3+C_4+C_5+C_6+C_7-C_8&=&\SI{-152}{\mV}\\
C_1+C_2+C_3+C_4+C_5+C_6-C_7+C_8&=&\SI{63}{\mV}\\
C_1+C_2+C_3+C_4+C_5-C_6-C_7-C_8&=&\SI{98}{\mV}\\
C_1+C_2+C_3+C_4-C_5-C_6+C_7+C_8&=&\SI{-35}{\mV}\\
C_1+C_2+C_3-C_4-C_5-C_6+C_7-C_8&=&\SI{0}{\mV}\\
C_1+C_2-C_3-C_4-C_5+C_6-C_7+C_8&=&\SI{145}{\mV}\\
C_1-C_2-C_3-C_4-C_5+C_6-C_7-C_8&=&\SI{-14}{\mV}
  \end{aligned}
  \qquad \right\}
  \tag{\ref{eq1}}
\end{equation}

In this system of equations (\ref{eq1}) there are only 8 coefficients: $C_i$, for $i=1,2,\ldots, 8$.

Hence this procedure makes it possible to compute, for each time series composed of a sequence of $n$ measured values, the values of the $n$ coefficients  $C_i$, for $i=1,2,\ldots, 8$, appearing in the system of $n$ linear algebraic equations that is generated when analyzing the time series by using the SWM. Equation (\ref{eq3}) makes it possible to compute the frequency of each of the $n$ trains of square waves $S_1,S_2,\ldots,S_n$. 
 
\section{Anlysis of a time series generated especially to illustrate essential features of the SWM}

The objective of this section is to emphasize the following characteristic of the SWM as a tool for the analysis of time series: The method remains applicable even when the changes in the value of the variable from which the sequence of samples is obtained are entirely unforeseeable; and the approximation provided for the measured values has a very high quality (i. e., the difference between any of the measured values and the corresponding computed value is very slight).

A random number generator (the MRNG \cite{b4,b5,b6}) was used to generate $10,000$ values between $-99.99999$ and $99.99999$.

When using the MRNG, the probability of any particular digit $0,1,2,\ldots,9$ appearing is the same: $\tfrac{1}{10}$.

To generate each one of those 10,000 numerical values, the following was done: If the first digit generated by the MRNG was between 0 and 4, inclusive, it was admitted that the numerical value generated (NVG) was negative. If, however, that digit was between 5 and 9, inclusive, it was accepted that the NVG was positive. The second digit of those generated was taken as the first digit of the whole part of the NVG. The third digit generated was taken as the second digit of the whole part of the NVG. The fourth, fifth, sixth, seventh and eighth digits generated were taken as the first, second, third, fourth and fifth digits, respectively, of the decimal part of the NVG. 

Each of the remaining values in the sequence of 10,000 numerical values to be analyzed with the SWM was obtained using the same type of procedure.

Of course, a different sequence of 8 digits was generated by the MRNG to be used to obtain each of those 10,000 numerical values.
 
It was supposed that the sequence of 10,000 values obtained by using the MRNG corresponded to a time series of values measured during a lapse of 5 seconds (\SI{5}{\second}) with a sampling frequency of \SI{2,000}{\hertz}. In other words, the MRNG was used to simulate the 10,000 consecutive measurements of a time series composed of a sequence of 10,000 numerical values. 

\begin{figure}[H]
\centering
\includegraphics[width=4in]{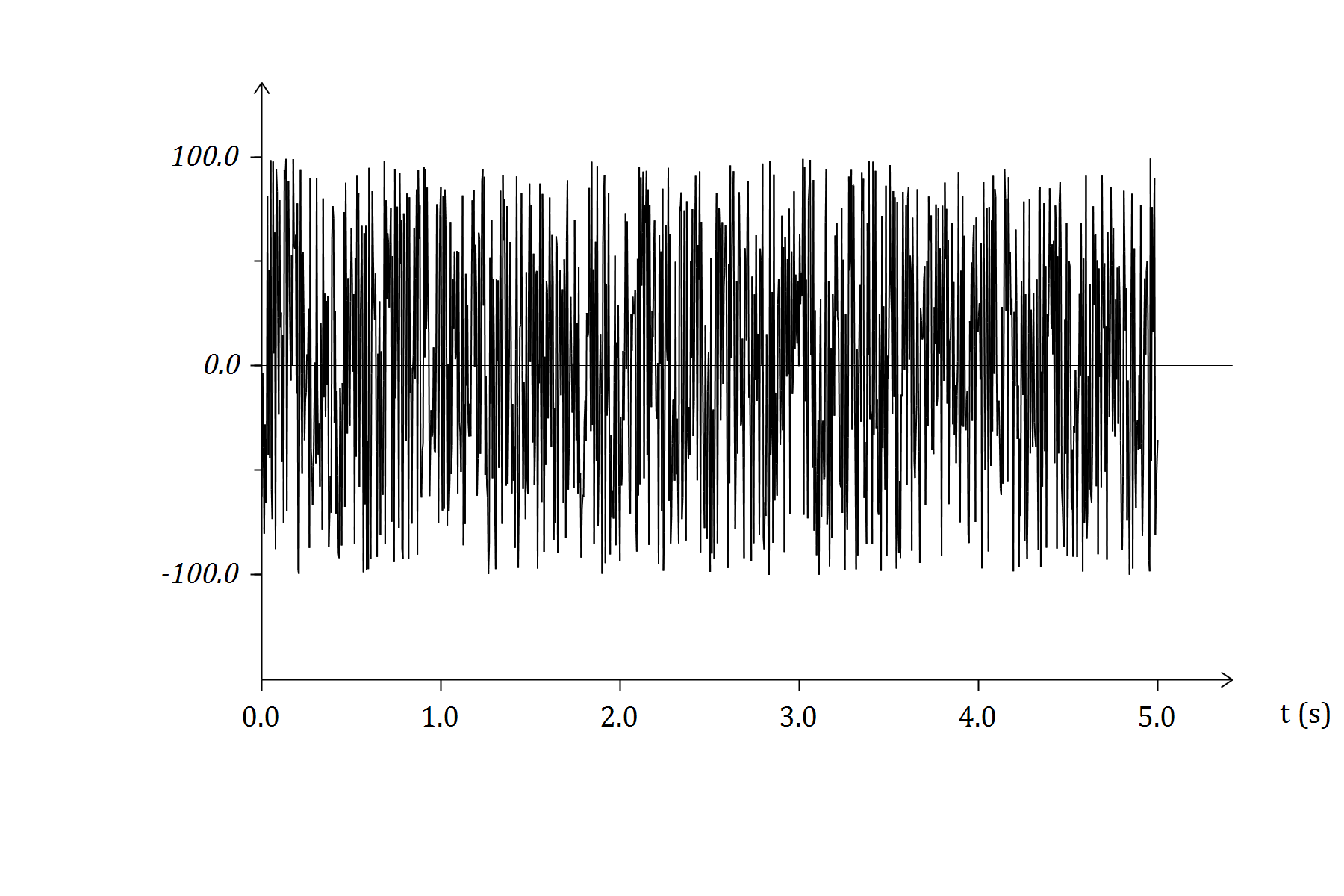}
\caption{The \textit{simulated} time series of 10,000 numerical values}
\label{f5}
\end{figure}

The \textit{simulated} time series described above is displayed in figure \ref{f5}.

The first 100 values of that numerical series of 10,000 are shown below:\\
\begingroup\makeatletter\def\f@size{8.5}\check@mathfonts
\begin{tabular}{r r r r r r r r}
$V_{1}=$& $\SI{-62.17387}{mV}$  & $V_{26}=$& $\SI{62.78057}{mV}$  & $V_{51}=$& $\SI{-58.04794}{mV}$  & $V_{76}=$& $\SI{-47.12626}{mV}$ \\ 
$V_{2}=$& $\SI{-77.8189}{mV}$  & $V_{27}=$& $\SI{52.57981}{mV}$  & $V_{52}=$& $\SI{-36.46843}{mV}$  & $V_{77}=$& $\SI{16.5785}{mV}$ \\ 
$V_{3}=$& $\SI{-86.98077}{mV}$  & $V_{28}=$& $\SI{-80.03298}{mV}$  & $V_{53}=$& $\SI{19.38409}{mV}$  & $V_{78}=$& $\SI{38.43836}{mV}$ \\ 
$V_{4}=$& $\SI{-30.27255}{mV}$  & $V_{29}=$& $\SI{-79.28799}{mV}$  & $V_{54}=$& $\SI{-37.07554}{mV}$  & $V_{79}=$& $\SI{-53.78808}{mV}$ \\ 
$V_{5}=$& $\SI{-21.89399}{mV}$  & $V_{30}=$& $\SI{-76.17741}{mV}$  & $V_{55}=$& $\SI{-32.84694}{mV}$  & $V_{80}=$& $\SI{-78.11808}{mV}$ \\ 
$V_{6}=$& $\SI{-65.28596}{mV}$  & $V_{31}=$& $\SI{36.22582}{mV}$  & $V_{56}=$& $\SI{-88.00845}{mV}$  & $V_{81}=$& $\SI{-80.16334}{mV}$ \\ 
$V_{7}=$& $\SI{-97.41397}{mV}$  & $V_{32}=$& $\SI{86.19701}{mV}$  & $V_{57}=$& $\SI{-77.15092}{mV}$  & $V_{82}=$& $\SI{46.22563}{mV}$ \\ 
$V_{8}=$& $\SI{70.16625}{mV}$  & $V_{33}=$& $\SI{-97.55277}{mV}$  & $V_{58}=$& $\SI{-73.22327}{mV}$  & $V_{83}=$& $\SI{16.87193}{mV}$ \\ 
$V_{9}=$& $\SI{-94.51527}{mV}$  & $V_{34}=$& $\SI{96.10761}{mV}$  & $V_{59}=$& $\SI{-86.89494}{mV}$  & $V_{84}=$& $\SI{22.48866}{mV}$ \\ 
$V_{10}=$& $\SI{-3.29921}{mV}$  & $V_{35}=$& $\SI{-28.53992}{mV}$  & $V_{60}=$& $\SI{15.32544}{mV}$  & $V_{85}=$& $\SI{81.87053}{mV}$ \\ 
$V_{11}=$& $\SI{27.59454}{mV}$  & $V_{36}=$& $\SI{-20.97347}{mV}$  & $V_{61}=$& $\SI{-17.86576}{mV}$  & $V_{86}=$& $\SI{-27.30575}{mV}$ \\ 
$V_{12}=$& $\SI{-76.21755}{mV}$  & $V_{37}=$& $\SI{-27.7922}{mV}$  & $V_{62}=$& $\SI{82.61901}{mV}$  & $V_{87}=$& $\SI{-33.18708}{mV}$ \\ 
$V_{13}=$& $\SI{-74.54863}{mV}$  & $V_{38}=$& $\SI{-34.08936}{mV}$  & $V_{63}=$& $\SI{49.45424}{mV}$  & $V_{88}=$& $\SI{39.52277}{mV}$ \\ 
$V_{14}=$& $\SI{-11.59495}{mV}$  & $V_{39}=$& $\SI{58.72335}{mV}$  & $V_{64}=$& $\SI{81.66214}{mV}$  & $V_{89}=$& $\SI{28.95455}{mV}$ \\ 
$V_{15}=$& $\SI{47.40731}{mV}$  & $V_{40}=$& $\SI{59.67262}{mV}$  & $V_{65}=$& $\SI{60.47497}{mV}$  & $V_{90}=$& $\SI{-40.53147}{mV}$ \\ 
$V_{16}=$& $\SI{5.70653}{mV}$  & $V_{41}=$& $\SI{20.80554}{mV}$  & $V_{66}=$& $\SI{41.39753}{mV}$  & $V_{91}=$& $\SI{-43.82255}{mV}$ \\ 
$V_{17}=$& $\SI{16.69057}{mV}$  & $V_{42}=$& $\SI{23.17323}{mV}$  & $V_{67}=$& $\SI{-78.66478}{mV}$  & $V_{92}=$& $\SI{-54.56893}{mV}$ \\ 
$V_{18}=$& $\SI{-83.79124}{mV}$  & $V_{43}=$& $\SI{9.25608}{mV}$  & $V_{68}=$& $\SI{-70.65624}{mV}$  & $V_{93}=$& $\SI{-47.19443}{mV}$ \\ 
$V_{19}=$& $\SI{-38.93865}{mV}$  & $V_{44}=$& $\SI{-39.10634}{mV}$  & $V_{69}=$& $\SI{-83.05503}{mV}$  & $V_{94}=$& $\SI{-56.6021}{mV}$ \\ 
$V_{20}=$& $\SI{42.47807}{mV}$  & $V_{45}=$& $\SI{40.03058}{mV}$  & $V_{70}=$& $\SI{-31.61922}{mV}$  & $V_{95}=$& $\SI{-22.46162}{mV}$ \\ 
$V_{21}=$& $\SI{76.04696}{mV}$  & $V_{46}=$& $\SI{-65.25953}{mV}$  & $V_{71}=$& $\SI{-66.50127}{mV}$  & $V_{96}=$& $\SI{38.18652}{mV}$ \\ 
$V_{22}=$& $\SI{-74.16542}{mV}$  & $V_{47}=$& $\SI{42.94391}{mV}$  & $V_{72}=$& $\SI{-49.22007}{mV}$  & $V_{97}=$& $\SI{73.71899}{mV}$ \\ 
$V_{23}=$& $\SI{99.33202}{mV}$  & $V_{48}=$& $\SI{-93.64582}{mV}$  & $V_{73}=$& $\SI{-42.51459}{mV}$  & $V_{98}=$& $\SI{27.59394}{mV}$ \\ 
$V_{24}=$& $\SI{86.24249}{mV}$  & $V_{49}=$& $\SI{3.09509}{mV}$  & $V_{74}=$& $\SI{1.38873}{mV}$  & $V_{99}=$& $\SI{0.41868}{mV}$ \\ 
$V_{25}=$& $\SI{-43.37761}{mV}$  & $V_{50}=$& $\SI{-1.77666}{mV}$  & $V_{75}=$& $\SI{-78.32747}{mV}$  & $V_{100}=$& $\SI{98.63304}{mV}$ \\ 
\end{tabular}
\endgroup
 \\
 \\

The first 100 dyads $(f_i,C_i)$, for $i=1,2,\ldots, 100$ (resulting from the analysis using the SWM) of the time series of 10,000 values is as follows:\\
\\
\begin{tabular}{ l l l l }
1. $(0.100000; -18065.729975)$   & 34.   $(0.100331; -48.628120)$   & 67.   $(0.100664; -30.526455)$ \\
2. $(0.100010; -85.317955)$   & 35.   $(0.100341; 93.884360)$   & 68.   $(0.100675; -209.547145)$ \\
3.  $(0.100020; -100.109675)$   & 36.   $(0.100351; 33.974930)$   & 69.   $(0.100685; -47.459485)$ \\
4.  $(0.100030; -91.989130)$   & 37.   $(0.100361; -289.039710)$   & 70.   $(0.100695; 31.143430)$ \\
5.  $(0.100040; 768.823705)$   & 38.   $(0.100371; -154.583500)$   & 71.   $(0.100705; -4.459135)$ \\
6.  $(0.100050; -54.968935)$   & 39.   $(0.100381; 3.959370)$   & 72.   $(0.100715; 2.026140)$ \\
7.  $(0.100060; 142.557690)$   & 40.   $(0.100392; 50.950590)$   & 73.   $(0.100725; -137.505715)$ \\
8.  $(0.100070; 40.323305)$   & 41.   $(0.100402; 1268.174905)$   & 74.   $(0.100735; 74.145900)$ \\
9.  $(0.100080; -184.475615)$   & 42.   $(0.100412; -47.335845)$   & 75.   $(0.100746; 119.423625)$ \\
10.  $(0.100090; -69.523010)$   & 43.   $(0.100422; 114.977245)$   & 76.   $(0.100756; 44.214730)$ \\
11.  $(0.100100; 325.761985)$   & 44.   $(0.100432; 24.150720)$   & 77.   $(0.100766; -3.829490)$ \\
12.  $(0.100110; 55.551825)$   & 45.   $(0.100442; 237.504795)$   & 78.   $(0.100776; -1.884085)$ \\
13.  $(0.100120; -134.091175)$   & 46.   $(0.100452; 28.773410)$   & 79.   $(0.100786; 6.228255)$ \\
14.  $(0.100130; 47.942450)$   & 47.   $(0.100462; 134.511410)$   & 80.   $(0.100796; 24.658120)$ \\
15.  $(0.100140; 27.621050)$   & 48.   $(0.100472; -0.563565)$   & 81.   $(0.100806; 290.612220)$ \\
16.  $(0.100150; 22.796750)$   & 49.   $(0.100482; -1150.327445)$   & 82.   $(0.100817; -75.545490)$ \\
17.  $(0.100160; 13883.376575)$   & 50.   $(0.100492; 22.801230)$   & 83.   $(0.100827; -51.233725)$ \\
18.  $(0.100170; -12.899225)$   & 51.   $(0.100503; 38.280880)$   & 84.   $(0.100837; -184.019420)$ \\
19.  $(0.100180; -328.522275)$   & 52.   $(0.100513; -20.427170)$   & 85.   $(0.100847; -188.399045)$ \\
20.  $(0.100190; 29.374290)$   & 53.   $(0.100523; -17.161430)$   & 86.   $(0.100857; -55.360550)$ \\
21.  $(0.100200; -108.606305)$   & 54.   $(0.100533; -52.007925)$   & 87.   $(0.100867; -7.625295)$ \\
22.  $(0.100210; -85.847715)$   & 55.   $(0.100543; -11.471110)$   & 88.   $(0.100878; -55.594890)$ \\
23.  $(0.100220; 187.059175)$   & 56.   $(0.100553; 280.304525)$   & 89.   $(0.100888; -1279.130445)$ \\
24.  $(0.100231; 12.679740)$   & 57.   $(0.100563; 458.847200)$   & 90.   $(0.100898; 124.454525)$ \\
25.  $(0.100241; -109.398105)$   & 58.   $(0.100573; 77.849500)$   & 91.   $(0.100908; -85.243345)$ \\
26.  $(0.100251; 108.646820)$   & 59.   $(0.100583; -33.323975)$   & 92.   $(0.100918; -145.174030)$ \\
27.  $(0.100261; -38.852280)$   & 60.   $(0.100594; -25.913470)$   & 93.   $(0.100929; 87.489570)$ \\
28.  $(0.100271; -83.227525)$   & 61.   $(0.100604; -523.193475)$   & 94.   $(0.100939; -94.767190)$ \\
29.  $(0.100281; 399.458855)$   & 62.   $(0.100614; 129.402450)$   & 95.   $(0.100949; 208.147540)$ \\
30.  $(0.100291; -156.149795)$   & 63.   $(0.100624; -98.541910)$   & 96.   $(0.100959; 15.027625)$ \\
31.  $(0.100301; -105.862070)$   & 64.   $(0.100634; 138.864805)$   & 97.   $(0.100969; -213.984295)$ \\
32.  $(0.100311; -22.212100)$   & 65.   $(0.100644; -1692.771570)$   & 98.   $(0.100980; 88.275250)$ \\
33.  $(0.100321; 910.844845)$   & 66.   $(0.100654; 3.439290)$   & 99.   $(0.100990; -39.121275)$ \\
                             &                                  & 100.      $(0.101000; -17.17885)$  \\

\end{tabular}
\\

\begin{figure}[H]
\centering
\includegraphics[width=4in]{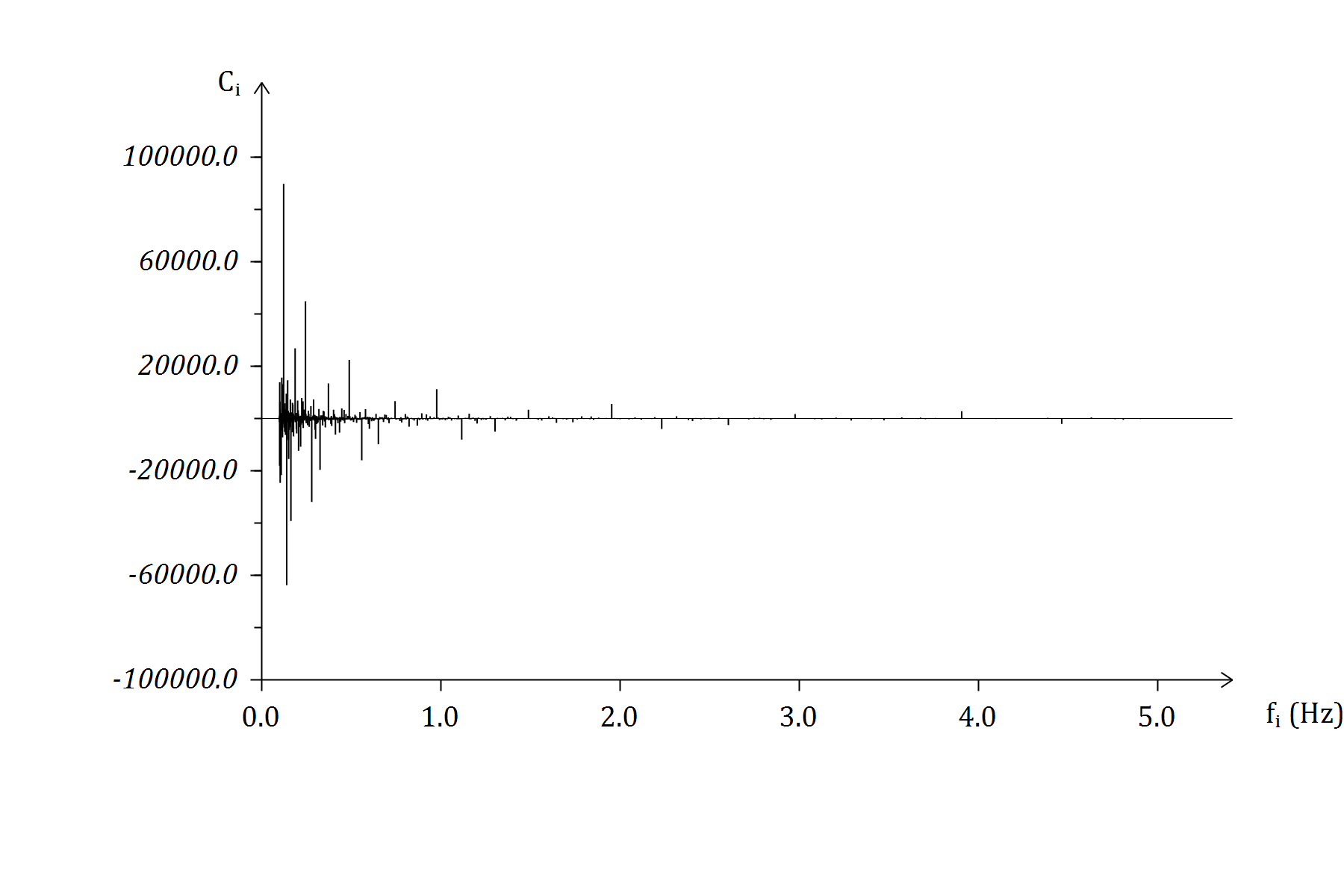}
\caption{Partial view of the SWT corresponding to the analysis (using the SWM) of the time series considered with 10,000 numerical values}
\label{f6}
\end{figure}

A partial display of the results obtained by using the SWM is presented in figure \ref{f6}.

The complete specification of the time series of 10,000 numerical values analyzed by the SWM and the corresponding sequence of 10,000 dyads of the type ($f_i$, $C_i$), for $i=1,2,3,\ldots, 10,000$, can be found in supplementary material available in the official page of the group \cite{h0}.

The approximation achieved (with the SWM) to each of the numerical values of the time series is outstanding. $V_i$, for $i=1,2,\ldots, 10,000$, will designate the $i$-th of those values, and $V_{i_{comp}}$ will refer to the corresponding approximation obtained by the SWM. $\abs{V_i-V_{i_{comp}}}$ is the absolute value of the difference between the two values. Of course, there will be 10,000 absolute values of this type. The maximum of those absolute values can be computed:
 \begin{equation*}
 \max \abs{V_i-V_{i_{comp}}} = 0.0000000000009379.
 \end{equation*}
\section{Discussion and prospects}
Consideration has been given to the application of the SWM to time series such that the value of the time interval between any two consecutive values remains constant. Specialists in a certain type of signal usually have the necessary tools to determine the minimum sampling frequency required to obtain a sequence of measured values of that type of signal, during a specific time interval, such that the sequence of values, or time series, can be considered to be an acceptable digital version of that signal, in that interval. 
Suppose that a particular time series is obtained by using that minimum sampling frequency or above. In that case, the analysis with the SWM of that time series is also a valid analysis of the corresponding signal, in that interval. 

Section 4 of this article provides strong support for the criterion that the SWM may be successfully applied even to acceptable digital versions of signals that have many abrupt changes.

The SWT source code in MATLAB is available on the official page of the group \cite{h1}, and in the official GitHub repository \cite{h2}. 

Future articles will deal with the following topics: 
\begin{enumerate}[I]
\item A general systematic approach for the use of the SWM for the analysis of functions of $n$ variables, for $n=1,2,3,\ldots$; and

\item Ways in which the use of the SWM can contribute to 1) the elimination of noise in different types of signals, and 2) the compression of information.
\end{enumerate}

\end{document}